\documentclass[11pt, reqno]{amsart}
\usepackage{amsmath, amssymb,amsfonts}
\usepackage[mathscr]{eucal}
\newtheorem{thm}{Theorem}[section]

\newtheorem{prop}[thm]{Proposition}
\newtheorem{lem}[thm]{Lemma}

\newtheorem{cor}[thm]{Corollary}

\newtheorem{conjecture}[thm]{Conjecture}

\def\Q{\mathbb Q}

\def\mod{\rm mod \, }

\def\mod{\operatorname{mod}}

\def\exp{\operatorname{exp}}

\def\li{\operatorname{li}}

\def\log{\operatorname{log}}

\title[A Generalised Kummer's Conjecture]
{A Generalised Kummer's Conjecture}

\date{\today}

\author{MJR Myers}

\address{Department of Mathematics and Statistics, 
Calvin College, Grand Rapids, Michigan  49546, USA}

\curraddr{School of Mathematical Sciences, University College Dublin, Belfield, Dublin 4, Ireland}

\email{mjm49@calvin.edu}

\thanks{The author would like to thank Robert Osburn and Larry Washington, and to give special thanks to M. Ram Murty. 
The author would also like to thank Pieter Moree at MPI
for his helpful comments.
This work was partially supported by IRCSET's Embark Postdoctoral Fellowship Scheme.}


\begin{document}

\begin{abstract}

Kummer's Conjecture predicts the rate of growth of the relative class numbers of cyclotomic fields of 
prime conductor. We extend Kummer's Conjecture to cyclotomic fields of conductor $n$, where $n$ is any natural number. We show that 
the Elliott-Halberstam Conjecture implies that this Generalised Kummer's Conjecture is true for almost all $n$ but is false for infinitely many $n$. 

\end{abstract}

\maketitle

\section{Introduction}

Let $\Q (\zeta_m) $ be the $m$th cyclotomic field, where $\zeta_m$ is a primitive $m$th root of unity for an integer $m \geq 1$. 
Let $h_m$ denote the class number of $\Q (\zeta_m)$ and 
$h_m^+$ be the class number of its maximal real subfield $\Q (\zeta_m + \zeta_m^{-1})$. 

Kummer proved that  the relative class number 
$h^-_{m} = h_{m} / h^+_{m} $ is an integer, and 
in 1851 he claimed (\cite{ku}, pg. 473) that the rule for the asymptotic growth of $h_p^-$ 
as the prime 
$p \to \infty$ is given by the formula
$$ 
\frac {p^{(p+3)/4}} {2^{(p-3)/2} \pi^{(p-1)/2}} =: G(p). \eqno (1)
$$
Kummer never published a proof of his claim, and the modern, rigourous reading of Kummer's assertion, that  
$$
\lim_{p \to \infty} \frac{h_p^-}{G(p)} =1
$$
has become well-known as ``Kummer's Conjecture''.

As it stands, Kummer's Conjecture remains unproven; however,  
Ankeny and Chowla \cite{ac} showed that 
$$\log (h^-_p / G(p)) = o(\log p )$$
as $p  \to \infty.$ 
Murty and Petridis  \cite{mp} proved what they called the Weak Kummer's 
Conjecture. They showed that  

there exists a positive constant $c$ such that 
$$c^{-1} \leq \frac {h^-_p} {G(p)} \leq c$$
holds for a sequence of primes $p_i$, where the number of primes $p_i \leq x$ is asymptotic 
to $x /  \log{x}$ as $x \to \infty$. 
With the additional assumption of the Elliott-Halberstam Conjecture, they were able to prove a stronger result. 
Recall that this conjecture says 
\begin{conjecture}[Elliott-Halberstam Conjecture] 
\label{eh}
For any $\delta >0$ and any $A>0$, 
$$
\sum_{k < x^{1 - \delta}} \max_{(l,k)=1} \max_{y \leq x} 
\left| \pi (y, k, l) - \frac{\li y}{\phi (k)} \right| \ll_{\delta, A}
\frac{x}{\log^A x}
$$
where $\pi (y, k, l)$ equals the number of primes $p \leq y$ such that $p \equiv l \mod  k $,
and $\li y = \int_2^y \frac{dt}{\log t}$. 
\end{conjecture}

Murty and Petridis showed that the Elliott-Halberstam Conjecture implies that for every $\epsilon > 0$ there exists an $x_{\epsilon}$
such that  
$$ 1 - \epsilon < \frac{h^-_p}{G(p)} < 1 + \epsilon $$
holds for all primes $x_{\epsilon} < p \leq x$,  with the exception of a set $P(\epsilon)$ such 
that 
$$| \{ p \in P(\epsilon): x_\epsilon < p \leq x \} | =o(\pi(x)) \text{.}$$

Hence Kummer's Conjecture concerning class numbers of cyclotomic fields is related to 
the density of primes in arithmetic progressions. Kummer's Conjecture is 
also related to pairs of primes. The Hardy-Littlewood Conjecture posits the existence of 
$\gg x / \log^2 x $ primes $p \leq x$ such that $2p+1$ is also prime; in 
1990 
Granville \cite{gr} 
proved that the Elliott-Halberstam Conjecture and the Hardy-Littlewood 
Conjecture together imply that Kummer's Conjecture is false. 
In that same paper Granville offered heuristic reasoning for believing that 
for all primes $p$ 
$$
(\log \log p )^{-1/2 +o(1)} \leq h_p^- /G(p) \leq  (\log \log p )^{1/2 +o(1)},
$$
and that these bounds are the best possible. 

More recently, Lu and Zhang \cite{lz} proved that for any fixed $\epsilon >0$, there is a positive number $Q$ 
depending only on $\epsilon$ 
such that for all primes $p \geq Q$, 
$$
e^{-1.4} p^{-\epsilon} (\log p)^{-1/3} \leq h_p^- / G(p) \leq e^{0.84} p^{\epsilon} (\log p)^{1/6} .
$$

In this paper we extend Kummer's Conjecture to composite numbers; that is, for natural numbers $n$ and a suitable 
function $G(n)$
(see (\ref{Gofn})), the Generalised Kummer's Conjecture predicts that $\lim_{n \to \infty} h_n^- / G(n) = 1$. 
We prove a composite moduli analogue of Murty and Petridis' Weak Kummer's Conjecture: 
\begin{thm}
\label{nnwkcuv}
Let $\omega (n)$ the number of distinct prime divisors of $n$. Then
$$
e^{- \omega (n)}  \ll \frac{h_n^-}{G(n)} \ll e^{ \omega (n)}
$$
holds for all but $o (x)$ natural numbers $n \leq x$.
\end{thm}

Moreover, assuming the Elliott-Halberstam Conjecture, the Generalised Kummer's Conjecture is true for 
almost all $n$ and is false for infinitely many $n$. More precisely we have the following two results.  
\begin{thm}
\label{nnwkccv}
Assume the Elliott-Halberstam Conjecture. Then for every $\epsilon > 0$ there exists an $x_{\epsilon}$
such that 
$$
1 - \epsilon < \frac{h_n^-}{G(n)} < 1 + \epsilon
$$
holds for all natural numbers $n \geq x_{\epsilon}$ with the exception of 
$ o(x) $ natural numbers $n<x$.
\end{thm}

\begin{thm}
\label{gkcf}
Assume the Elliott-Halberstam Conjecture. 
Then the Generalized Kummer's Conjecture fails 
for infinitely many natural numbers $n$. 
\end{thm}

\section{Generalised Kummer's Conjecture}

For the cyclotomic field $\Q(\zeta_n)$, 
one may obtain the formula 
(see pg. 42 of \cite{wa}): 
\begin{equation}
\label{poleratio}
h_n^- = \frac{Q w \sqrt{|d_n / d_n^+|}}{\pi^{\phi (n) /2} 2^{ \phi (n)/2}} 
\prod_{\substack{\chi \ \mod n \\ \chi \ \text{odd} }} L(1, \chi ).
\end{equation}
Here $d_n$ is the discriminant of $\Q(\zeta_n)$ and $d_n^+$ is the discriminant of $\Q (\zeta_n + \zeta_n^{-1})$. 
Also, $w$ is the number of roots of unity in $\Q(\zeta_n)$, and 
$Q=1$ if $n$ is a prime power $p^r$ and $Q= 2$ otherwise.  

By   Proposition 2.7 in  \cite{wa}, 
\begin{equation*} \label{wa 2.7}
d_n  = (-1)^{\phi (n) /2} \frac{n^{\phi (n)}}{\prod_{p|n} p^{\phi (n) / (p-1)}},
\end{equation*}
and by Lemma 4.19 in \cite{wa}, 
\begin{equation*}
\label{wa 4.19}
d_n  = \begin{cases}
\ \  p(d_n^+)^2  & \text{if $n = p^r$ with  $p \ne 2$},\\
\ \  4(d_n^+)^2  & \text{if $n = 2^r$},\\
\ \  (d_n^+)^2    & \text{otherwise.}
\end{cases}
\end{equation*}

Hence we see that 
\begin{equation*}
h_n^- = 
a_n \left( \frac{1}{2 \pi} \sqrt{ \frac{n}{\prod_{p | n}p}} \right)^{\phi(n) /2}
\prod_{\substack{\chi \ \text{odd} \\ \chi \mod n }} L(1, \chi ),
\end{equation*}
with \begin{equation*}
a_n  = \begin{cases}
\ \  2p^{r +1/4}  & \text{if $n = p^r$ with  $p \ne 2$},\\
\ \  2^{r+ 1/2}  & \text{if $n = 2^r$},\\
\ \  4n & \text{if odd $n \ne p^r$ }, \\
\ \  2n  & \text{if even $n\ne 2^r$ .}
\end{cases}
\end{equation*}
Let 
\begin{equation}
\label{Gofn}
G(n) = a_n \left( \frac{1}{2 \pi} \sqrt{ \frac{n}{\prod_{p | n} p}} \right)^{\phi(n) /2}.
\end{equation}
Then the composite moduli form of Kummer's Conjecture may be stated as follows.
\begin{conjecture}[Generalized Kummer's Conjecture]
\label{gkc} 
$$h_n^- \sim G(n)$$ as the natural number $n \to \infty$. 
\end{conjecture}

As in \cite{mjrm1}, we can rewrite the product of $L$-functions in (\ref{poleratio}) as 
$$\prod_{\substack{\chi \ \mod n \\ \chi \ \text{odd} }} L(1, \chi ) =
\exp \left( \frac{\phi (n)}{2} f_n \right), 
$$
where $f_n = \lim_{x \to \infty} f_n (x)$ 
and $f_n (x)$ is the finite sum 
$$f_n (x) = \sum_{r \leq x} \frac{c_n (r)}{r}$$ 
with 
$$ c_{n} (r) = \begin{cases}
\ \  1 & \text{if $r = q^m \equiv 1 \mod n$},\\
-1 & \text{if $r= q^m \equiv -1 \mod n$},\\
\ \  0 & \text{otherwise}
\end{cases} $$
where $q$ is a prime and $m \geq 1$.

Clearly 
the Generalized Kummer's Conjecture is true if and only if $f_n = o\left( \frac{1}{\phi (n)} \right).$

\section{Lemmas}

We will need the following theorems which are also used in \cite{mp}.

\begin{lem}[Siegel-Walfisz Theorem]
\label{sw}
For any constant $A>0$, there is a constant $c(A) > 0$
so that uniformly for $x \geq 3$, $1 \leq k \leq (\log x)^A, 
(k, l)=1$, we have 
$$\pi (x, k, l) = \frac{ \li x}{ \phi (k)} + O \left( x e^{-c(A) \sqrt{\log x}} \right).$$
\end{lem}

\begin{lem}[Bombieri-Vinogradov Theorem]
\label{bv}
Assume $x \geq 2$.  
For any $A > 0$ there exists $B= B(A) > 0$ such that 
$$ \sum_{k < \frac{x^{1/2}}{\log^B x }} E(x, k) \ll_A \frac{x}{\log^A x} $$
where 
$$E(x, k) := \max_{y \leq x} \max_{(k,l) = 1} \left| \pi (y, k, l) - \frac{\li y}{\phi (k)} \right|.
$$
\end{lem}

\begin{lem}[Brun-Titchmarsh Theorem]
\label{bt}
For $k<x$, $(k,l) = 1$, 
$$ \pi (x, k, l) < \frac{2x}{\phi (k) \log (x/k)} .$$
\end{lem}

\begin{lem}[\cite{ho}, pg. 124]
\label{hooley}
Let $l$ be a fixed, non-zero integer, and $ \epsilon , A, B$ positive real numbers 
where $A > B +30$. 
Then for any numbers $x$ and $X$ such that $x^{1/2} < X < x \log^{-A} x$ and 
$x> x_0 (\epsilon, B)$, we have
$$ \pi (x, k, l) \leq \frac{(4 + \epsilon )x}{\phi (k) \log X}$$
for every $k$ such that $X \leq k \leq 2X $, and $(l, k)=1$, except for at most 
$X \log^{-B} x$ exceptional values of $k$.
\end{lem}

\begin{lem}[\cite{mp}, pg. 298]
\label{mp 2.3}
Fix $l$ and $k$, $(l, k) =1$. The number of primes $x < p \le 2x$ such that $kp+l$ 
is also prime is
$$
\ll \prod_{p_i | kl} \left( 1- \frac{1}{p_i} \right)^{-1} \frac{x}{\log^2 x}
$$
uniformly for $k<x^2$.
\end{lem}

\begin{lem}[\cite{mp}, pg. 298]
\label{mp 2.4}
There is a constant $c$ such that, as $T \to \infty$,
$$
\sum_{k \leq T} \prod_{p_i |k} \left( 1- \frac{1}{p_i} \right)^{-1} \sim cT.
$$
\end{lem}

\begin{lem} [\cite{mjrm1}, Corollary 3.6]
\label{murty cong cor}
The number of solutions $\mod n$ to $ x^m \equiv 1 \mod n$ (or to $ x^m \equiv -1 \mod n$) is at most $2m^{\omega (n)}$,
where $\omega (n) $ is the number of distinct prime divisors of $n$. 
\end{lem}

Let $A$ be some constant greater than or equal to $e$. A slight modification of Hardy and Ramanujan's original proof in \cite{hr2} of the normal order of $\omega (n)$  
shows that the number of $n \leq x$ such that $\omega (n) > A \log \log x$ is $o\left(\frac{x}{\log x}\right)$.
More specifically, we have 
\begin{lem}
\label{hr lemma}
For any constant $A \geq e$, 
$$|\{ n \leq x :  \omega (n) > A \log \log x \}|  \ll \frac{x}{(\log x)^{1 + A \log A - A}(\log \log x)^{1/2}}. $$
\end{lem}

\begin{proof}
Let  
$$
\pi (x, k) = \sum_{\substack {n \leq x \\ \omega(n) = k} } 1 \ \ \text{ and } \ \ S = \sum_{k \geq A \log \log x} \pi (x, k+1).
$$ 
Lemma B of \cite{hr2} gives us the uniform upper bound
$$
\pi (x, k+1) < \frac{Lx}{\log x} \ \frac{(\log \log x +D)^{k}}{k!},
$$
where $L$ and $D$ are absolute constants. Clearly, then   
$$
S <  \frac{Lx}{\log x} \ \sum_{k \geq A \log \log x} \frac{(\log \log x +D)^{k}}{k!}.
$$

Write $\xi = \log \log x +D$ and let $k_1$ be the smallest integer greater than $A\xi$. 
Then 
\begin{eqnarray*}
 \sum_{k \geq A \xi} \frac{\xi^k}{k!}
&<&  \ \frac{\xi^{k_1}}{k_1!} \left[ 1+\frac{\xi}{k_1 + 1} +\frac{\xi^2}{ (k_1+1)(k_1+2)} + \cdots \right] \\
&<& \frac{\xi^{k_1}}{k_1!} \left[ 1+\frac{1}{A} + \frac{1}{ A^2 } + \cdots \right]  
      = \ \frac{ \xi^{k_1}}{k_1!} \left[ \frac{A}{A-1} \right] \\
&\ll&  \frac{e^{k_1(\log \xi - \log k_1 + 1)}}{\sqrt{k_1}} 
\end{eqnarray*}
by Stirling's formula. It follows, then, that 
$$
\sum_{k \geq A \xi} \frac{\xi^k}{k!} \
 \ll \  \frac{e^{(A - A\log A)\xi}}{\sqrt{\xi}} \ \ll \ (\log x)^{A-A\log A} (\log \log x)^{-1/2}.
$$
Thus the number of $n \leq x$ such that $\omega (n) > A \log \log x$ is 
$$
 \ll \frac{x}{(\log x)^{1 + A \log A - A}(\log \log x)^{1/2}},
$$
as required.
\end{proof}

\section{Unconditional Composite Moduli Weak Kummer's Conjecture}

In this section we will prove Theorem \ref{nnwkcuv}, the 
Weak Kummer's Conjecture for composite moduli.
To remove the contributions of the prime powers $q^m$ with $m \geq 2$ from the sum $f_n$, 
we will use the following lemma:

\begin{lem}
\label{chuck powers 2}
$$\sum_{q  \text{ prime}} \sum_{m \geq 2 }   \frac{c_{n} (q^m)}{mq^m}
= O\left( \frac{\omega (n)}{n}   \right).$$
\end{lem}

\begin{proof}
Write $ S = \sum_{q  \text{ prime}} \sum_{m \geq 2 } \ \  \frac{c_{n} (q^m)}{mq^m}$. 
Then clearly 
$$| S |  \leq  \sum_{q  \text{ prime} } \sum_{\substack{ m \geq 2 \\ q^m \equiv \pm 1 (n) }}  \frac{1}{mq^m}  
= S_1 + S_2$$
where 
\begin{eqnarray*}
S_1 =   \sum_{\substack {q  \text{ prime} \\ q < n }}   \sum_{\substack{ m \geq 2 \\ q^m \equiv \pm 1 (n) }}   \frac{1}{mq^m}   
& \text{ and } &
S_2 = \sum_{\substack {q  \text{ prime} \\ q > n }}   \sum_{\substack{ m \geq 2 \\ q^m \equiv \pm 1 (n) }}   \frac{1}{mq^m}.
\end{eqnarray*}

Recall that $\pi (x)$ is the number of primes $p \leq x$. By the prime number theorem, we have 
\begin{eqnarray*}
|S_2| & \leq & \sum_{\substack {q > n \\ q \ \text{prime}}}  \sum_{m \geq 2 } q^{-m}  
              \ll   \sum_{\substack {q > n \\ q \ \text{prime}}} q^{-2}  \\   
& \ll &  \lim_{x \to \infty} \left[  \frac{\pi (t)}{t^2}  \right]_{n}^{x}  + \int_n^x \frac{\pi(t)}{t^3}dt     \\              
& \ll &  \frac{1}{n \log n} + \int_n^{\infty} \frac{1}{t^2 \log t}dt   \\
& \ll &  \frac{1}{n \log n} = o \left( \frac{1}{n} \right) .
\end{eqnarray*}
We now consider $S_1 = S_3 + S_4$ with 
\begin{eqnarray*}
S_3 = \sum_{m \geq 2 } \sum_{\substack{q < n \\ q^m \equiv 1 (n) }}   \frac{1}{mq^m}
 & \text{and} &
S_4 = \sum_{m \geq 2 } \sum_{\substack{q < n \\ q^m \equiv -1 (n) } }   \frac{1}{mq^m}.
\end{eqnarray*}

For a fixed $n$, let $C(m)$ denote the number of solutions $x<n$ to the congruence $x^m \equiv 1 \mod n$. 
By 
Lemma \ref{murty cong cor}  $ C(m) \leq 2m^{\omega (n)}$, where 
$\omega (n)$ is the number of distinct prime divisors of $n$. 
This gives us the upper bound
$$
\left| \bigcup_{i=2}^m \{x<n : x^i \equiv 1 \mod n \}\right| \leq \sum_{i=2}^m C(i) \leq \sum_{i=2}^m 2i^{\omega (n)} =:B(m).
$$
Observe that $B(m) \ll m^{\omega (n) +1}$.

Now for each solution $x<n$ to $x^m \equiv 1 \mod n$, write $x^m = u_i n+1$, where each $u_i$ is 
a distinct positive integer for $i = A(m), \dots, A(m) + C(m) - 1$. Here $A(2) = 1$, and for $m \geq 3$,
$A(m) = B(m-1) + 1$. Then  
$$S_3 = \sum_{m \geq 2 } \frac{1}{m} \sum_{\substack{q < n \\ q^m \equiv 1 (n)}}   \frac{1}{q^m} 
= \sum_{m \geq 2} \frac{1}{m} \left( \sum_{a= A(m)}^{B(m)} \frac{\theta_a}{ u_a n +1} \right),$$
where for $ A(m) \leq  a \leq A(m) +C(m) - 1$, 
\begin{equation*}
\theta_a  = \begin{cases}
\ \  1  & \text{if $\sqrt[m]{u_a n + 1} $ is a prime},\\
\ \  0    & \text{otherwise.}
\end{cases}
\end{equation*}
For any $ A(m) + C(m) \leq a \leq B(m)$, let $\theta_a = 0$ and $ u_a = 1$.
Thus 
$$ S_3  <  \frac{1}{n} \sum_{m \geq 2} \frac{1}{m} 
\left( \sum_{a= A(m)}^{B(m)} \frac{\theta _a}{ u_a} \right). 
$$

Let $d_m$ be the inner sum $ \sum_{a= A(m)}^{B(m)} \frac{\theta _a}{ u_a}$, and write the partial sums
$$D_2 = d_2, D_3 = d_2 + d_3, D_4 = d_2 + d_3 + d_4, \dots .$$ 
Notice that for any two indices $a \not= a'$ such that $\theta_a$ and $\theta_{a'}$ are both nonzero, 
we must have $u_a \not= u_{a'}$, which implies that 
$D_r \leq \sum_{a=1}^{B(r)} a^{-1}$.

Now,   
$$
\sum_{m=2}^x \frac{d_m}{m} = \frac{D_x}{x} + \sum_{r=2}^{x-1} \frac{D_r}{r(r+1)},
$$
and we get
\begin{eqnarray*}
0  <  S_3 
& = &  \frac{1}{n} \sum_{m \geq 2} \frac{D_m}{m(m+1)} \\
& \ll & \frac{1}{n} \sum_{m \geq 2} \frac{1}{m(m+1)} \left( \sum_{a=1}^{m^{\omega (n) +1}} \frac{1}{a} \right) \\
& \ll & \frac{\omega (n)}{n} \sum_{m \geq 2} \frac{\log (m)}{m(m+1)} 
 =  O \left( \frac{\omega (n)}{n} \right).
\end{eqnarray*}

Similarly we can express the sum 
$S_4 = \sum_{m \geq 2} \frac{1}{m} \left( \sum_{a= A(m)}^{B(m)} \frac{\phi_a}{ v_a n -1} \right)$
where the $v_a$'s are positive integers and 
$\phi_a = 1$ or $0$.
Since 
$$\sum_{a= A(m)}^{B(m)} \frac{\phi_a}{ v_a n -1} = \sum_{a= A(m)}^{B(m)} \frac{\phi_a}{ v_a n +1}
+ O(n^{-2}),$$ 
it follows that 
$S_4 \ll \frac{\omega (n)}{n}$ as well. 
\end{proof}

Define the sum
$$g_{n}(x) = \sum_{ \substack {q \ \text{prime } \\  q \leq x }} \frac{c_{n} (q)}{q},$$
and 
$$ 
g_{n} = \lim_{x \to \infty} g_{n} (x).
$$
By Lemma \ref{chuck powers 2}, $f_n = g_n + O(\frac{\omega (n)}{n})$. 
An application of the Siegel-Walfisz Theorem (Lemma \ref{sw}) reduces the infinite sum 
$g_n$ to a finite one. 

\begin{lem} 
$$g_{n} = g_{n} (2^{n}) + O (n^{-2}).$$ 
\label{sw lemma 2}
\end{lem}

\begin{proof}
For $x \geq y \geq 3$, Riemann-Stieltjes integration gives 
\begin{equation}
\label{rsi}
g_{n} (x) -g_{n} (y) = \sum_{y < q \leq x} \frac{c_{n} (q)}{q} 
= \left[ \frac{A_{n} (t)}{t} \right]_y^x + \int_y^x \frac{A_{n} (t)}{t^2} dt,
\end{equation}
where $A_{n} (t) = \pi (t, n, 1) - \pi (t, n, -1) $.
Using the Siegel-Walfisz Theorem and taking $x > 2^{n}$, we obtain 
$A_{n} (x)  \ll \frac{x}{n \log^2 x}
$
and so
\begin{eqnarray*}
| g_{n}(x) - g_{n}(2^{n}) | 
& \ll & \left[ \frac{A_{n} (t)}{t} \right]_{2^{n}}^x + \int_{2^{n}}^x \frac{A_{n} (t)}{t^2} dt  \\
& \ll & \frac{1}{n \log^2 x} - \frac{1}{n (\log 2^{n})^2} 
+ \frac{1}{n} \left[ \frac{1}{\log t} \right]_{2^{n}}^x \\
& \ll & \frac{1}{n^{2}}.
\end{eqnarray*}
\end{proof}

We have shown 
$$
f_n = g_n +O \left( \frac{\omega (n)}{n} \right) = g_n (2^n) + O \left( \frac{\omega (n)}{n} \right)
$$
and have reduced the problem to one of studying the finite sum 
$$g_{n}(2^n) = \sum_{\substack {q \ \text{prime, } \\ q \leq 2^n }} \frac{c_{n} (q)}{q}.$$
We will now find bounds on $g_n (2^n)$ by using (\ref{rsi}) to partition this sum into terms  
on which we may apply our various estimates for $\pi (t, n, 1) - \pi (t, n, -1) $. 
We are now in a position to prove Theorem \ref{nnwkcuv}. 

\begin{proof}

Note that 
$$
| A_{n} (t)|= |\pi (t, n, 1) - \pi (t, n, -1) | \leq 2 E(t, n)
$$
as defined in  Lemma \ref{bv},  
the conditions of which are satisfied 
for  $x < n \leq 2x$ and $n^{2} \log ^{2B} n < q < 2^{n}$.
Hence, 
\begin{eqnarray*}
\sum_{x < n \leq 2x} \sum_{ n^2 \log^{2B} n < q < 2^{n}} \frac{c_{n} (q)}{q}
& \ll & \left[ \frac{\sum_{x < n \leq 2x} E(t, n)}{t} \right]_{n^2 \log^{2B}n}^{2^{n}} 
+ \int_{n^2 \log^{2B}n}^{2^{n}} \frac{\sum_{x < n \leq 2x} E(t, n)}{t^2} dt \\
& \ll & \left[ \frac{1}{\log^A t} \right]_{x^2 \log^{2B} x}^{\infty} 
+ \int_{x^2 \log^{2B} x}^{\infty} \frac{1}{t \log^A t} dt \\
&\ll &  \log^{-A+1} x.
\end{eqnarray*}
If we set $D(n) =  g_{n} (2^{n}) - g_{n} (n^2 \log^{2B} n)$, then we have shown 
$$\sum_{x< n < 2x} |D(n) | \ll \frac{1}{\log^{A-1} x}.$$
Thus for any constant $c > 0$, 
\begin{eqnarray*}
\# \left\{ x < n \leq 2x \ : \ |D(n)| > \frac{c}{n} \right\} 
& \ll & \frac{x}{\log^{A-1} x}.
\end{eqnarray*}
Take $A > 3$ in Lemma \ref{bv}.

By dyadic decomposition we discard at most $x \log^{-A+1} x$ natural numbers 
and we now restrict our 
attention to primes $q$ in the range $ n^2/4 < q  \leq n^{2} \log^{2B} n $. 
By Lemma \ref{bt},  
\begin{eqnarray*}
\sum_{ n^{2}/4 < q \leq n^{2} \log^{2B} n } \frac{c_{n} (q)}{q} 
& \ll &  \frac{1}{\phi (n) \log n} + 
         \int_{n^{2}/4}^{n^{2} \log^{2B} n } \frac{1}{t \phi(n) \log (t/ n) } dt \\
& \ll & 
          \frac{1}{\phi(n)} \log \left( 1 + \frac{\log ( \log^{2B} n )}{\log (n /4)} \right) \\
& = & o \left(\frac{1}{\phi (n)} \right).
\end{eqnarray*}
That is, $g_{n} (n^{2} \log^{2B} n) - g_{n} (n^{2}/4) = o(1/n)$, and we 
now consider the range $ 2^A n \log^A n < q \leq  n^{2}/4$.

Take $X < n < 2X$ and  let
$2^A n \log^A n < t < n^{2}/4$.
The conditions of Lemma \ref{hooley} are satisfied in this range, and so 
\begin{eqnarray*}
\sum_{2^A n \log^A n < q \leq n^{2}/4 } \frac{c_{n} (q)}{q} 
&\ll & \frac{1}{\phi (n) \log n} + \frac{1}{\phi (n) } \int_{2^A n \log^A n }^{ n^{2}/4} \frac{1}{t \log t} dt \\
&  \ll &  \frac{1}{\phi (n)}.
\end{eqnarray*}
This holds for all natural numbers $X < n < 2X$ with the exception of a set of size
$\ll \frac{X}{\log^B X}$. Using dyadic decomposition, we see the number of
exceptional $n<x$ is $ \ll \frac{x}{\log^B x}$.

To estimate $g_n(2^A n \log^A n) - g_n ( n \log n) $, we apply Lemma \ref{bt} again to get 
\begin{eqnarray*} 
\sum_{ n \log n < q \leq 2^A n \log^A n  } \frac{c_{n} (q)}{q}
& \ll & \frac{1}{\phi (n)} \left[ \log \left( \log \frac{t}{n} \right) \right]_{n \log n}^{2^A n \log^A n} \\
& \ll & 
          \frac{1}{\phi (n)} \log \left( 1 + \frac{\log (2^A \log^{A-1} n) }{\log (\log n ) } \right) \\
& \ll &  \frac{1}{\phi (n)} .
\end{eqnarray*}

Using Lemma \ref{bt} one more time: 

\begin{eqnarray*}
\sum_{\epsilon n \log n / \log \log n < q < n  \log n } \frac{c_n(q)}{q}
& \ll & 
\int_{\epsilon n \log n / \log \log n}^{n  \log n} 
\frac{dt}{ \phi (n) t \log (t/n)} \\
& \ll & \frac{1}{\phi (n)} \log \left(   \frac{\log (\log n)}{\log ( \log n / \log \log n ) }      \right) 
 =  o \left( \frac{1}{\phi (n)}   \right).
\end{eqnarray*}

Finally, we need to analyse the sum
$$
g_n ( \epsilon n \log n / \log \log n ) = \sum_{n < q < \epsilon n \log n / \log \log n} \frac{c_n (q)}{q}.
$$
This sum is $\neq 0$ when there are summands; that is, when at least one of 
$
n \pm 1, 2n \pm 1, \dots , kn \pm 1 
$
is prime for $k < \epsilon \log x / \log \log x $. Hence we use the Prime Number Theorem for arithmetic progressions and we see 
\begin{eqnarray*}
\sum_{k < \frac{\epsilon \log x}{ \log \log x}}  \frac{k}{\phi (k)} \frac{x}{\log x}
& \ll & \frac{x}{\log x} \cdot \frac{\epsilon  \log x}{ \log \log x} 
 \ll  \frac{\epsilon x}{ \log \log x} = o(x).
\end{eqnarray*}
The number of $ n \leq x$ such that $g_n ( \epsilon n \log n / \log \log n ) \neq 0$ 
is $o(x)$, and Theorem \ref{nnwkcuv} has been proved.
\end{proof}

\section{Conditional Composite Moduli Weak Kummer's Conjecture}

We now prove Theorem \ref{nnwkccv} that Kummer's Conjecture holds for almost all $n$.  
We will need natural numbers analogues of 
Propositions 1 and 2 from 
\cite{gr}.

\begin{prop}
\label{grprop1nn}
$$
\sum_{m \geq 2} \frac{1}{m} \sum_{\substack{q^m \equiv \pm 1 \mod n \\ q \text{ prime} }} \frac{1}{q^m} = o \left( \frac{1}{n} \right)
$$
for all but $o\left( \frac{x}{\log x} \right)$ natural numbers $n \leq x$.
\end{prop}

\begin{proof}
For any prime $q> n$, 
\begin{eqnarray*} 
\sum_{m \geq 2} \frac{1}{mq^m} 
& \leq & \frac{1}{2} \left( \frac{1}{q^2} + \frac{1}{q^3} + \cdots \right) 
\leq \frac{1}{q^2}.
\end{eqnarray*}
Also, for any prime $q < n$,
$$
\sum_{m \geq 2} \frac{1}{mq^m} 
\leq \frac{1}{2n^{2}} \left( 1+ \frac{1}{q} + \frac{1}{q^2} + \cdots \right)
\leq \frac{1}{n^{2}}. 
$$
Thus if we list the primes $q_i$ in order so that $ q_k < n < q_{k+1}$,
\begin{eqnarray*} 
\sum_{m \geq 2} \frac{1}{m} \sum_{\substack{q^m \equiv \pm1 \mod n \\ q^m > n^{2} }} \frac{1}{q^m} 
& \leq & \sum_{m \geq 2} \frac{1}{m q_1^m} + \cdots + \sum_{m \geq 2} \frac{1}{m q_k^m} 
+ \sum_{m \geq 2} \frac{1}{m q_{k+1}^m} + \sum_{m \geq 2} \frac{1}{m q_{k+2}^m} + \cdots  \\
& \leq & \frac{1}{n^2} + \cdots \frac{1}{n^2}
+ \frac{1}{q_{k+1}^2} + \frac{1}{q_{k+2}^2} + \cdots \\
&=& O \left( \frac{1}{n^2} \cdot \frac{n}{\log n} \right)  = O \left( \frac{1}{n \log n} \right) 
\end{eqnarray*}
by the Prime Number Theorem.

Again using  Lemma \ref{murty cong cor}, we have $\leq 4(2)^{\omega (n)}$ solutions mod $n$ of the congruence
$x^m \equiv \pm 1 \mod n$. 
Also, by  Lemma \ref{hr lemma},
the number of $n \leq x$ such that $\omega (n) > A \log \log x$ is $o\left(\frac{x}{\log x}\right)$. (Here we choose 
$ e \leq A < 2 / \log 2$.)

Hence 
\begin{equation}
\label{damn it}
\frac{1}{2}\sum_{\substack{q^2 \equiv \pm1 \mod n \\ n \log^2 n < q^2 \leq n^2 }} \frac{1}{q^2} \ 
\leq \frac{1}{2n\log^2 n}   \sum_{\substack{q^4 \equiv 1 \mod n \\ q \leq n }} 1 \
\leq \ \frac{2(2^{\omega (n)})}{n \log^2 n} \ = \ o \left( \frac{1}{n} \right)
\end{equation}
for all but $o\left(\frac{x}{\log x}\right)$ numbers  $n \leq x$.

Now if $q^m \leq n^2$, then $m < 4 \log n$.
Let 
$$
S_{n} = \sum_{3 \leq m \leq 4\log n} \frac{1}{m} 
\sum_{\substack{q^m \equiv \pm 1 \mod n \\ q^m \leq n \log^2 n} } \frac{1}{q^m}. 
$$
Then 
\begin{eqnarray*}
\sum_{x <n <2x} S_{n} & \ll & \sum_{3 < m < \log n} \frac{1}{m} \sum_{x \log^2 x < q^m < x^2 } \frac{1}{q^m} 
                              \sum_{ \substack{ x<n<2x \\ n| q^m \pm 1 }} 1 \\
& \ll & \sum_{3 < m < \log n} \frac{1}{m} \sum_{\substack{q < x^{2/m} \\ q \text{ prime} }} \frac{1}{x \log^2 x} \cdot x^{\delta} \\
& \ll & \log \log x \cdot \frac{x^{2/3}}{\log (x^{2/3})} \cdot \frac{1}{x^{1-\delta} \log^2 x} \\
& \ll & \frac{\log \log x}{\sqrt[4]{x} \log^3 x},
\end{eqnarray*}
because the number of divisors of $q^m + 1$ or of $q^m - 1$ is $o(x^{\delta})$ for any $\delta >0$ (see \cite{ap}, pg. 296).

The number of $ x < S_n < 2x$ such that $S_n > \epsilon /n$ is $\ll x^{3/4}$, since otherwise
$$
\sum_{x <n <2x} S_{n} \gg \frac{1}{x} \cdot x^{3/4} = \frac{1}{\sqrt[4]{x}}. 
$$
By dyadic decomposition together with (\ref{damn it}), we see that
$$ 
\sum_{m \geq 2} \frac{1}{m} 
\sum_{\substack{q^m \equiv \pm 1 \mod n \\ n \log^2 n \leq q^m \leq n^2 } } \frac{1}{q^m} = o (1/n)
$$
for all but $o(x/\log x) $ natural numbers $n \leq x$.
It now suffices to show that
$$
s_{n} = \sum_{m \geq 2} \frac{1}{m} 
\sum_{\substack{q^m \equiv \pm 1 \mod n \\  q^m < n \log^2 n  } } \frac{1}{q^m} = o (1/n)
$$
for all but $o(x/\log x)$ numbers $n \leq x$.

Note that 

\begin{eqnarray*}
\sum_{x <n <2x} s_{n} & \ll & \sum_{m \geq 2} \frac{1}{m} \sum_{x < q^m < x \log^2 x} \frac{1}{q^m} 
                              \sum_{ \substack{ x<n<2x \\ q^{m}= \pm 1 + kn }} 1 \\
& \ll & \sum_{m \geq 2} \frac{1}{m} \sum_{x < q^m < x \log^2 x} \frac{\log^2 x}{q^m} \\
& \ll & \sum_{\substack{q \ \text{prime} \\ q < x^{1/2} \log x }} \frac{\log^2 x}{x} \\
& \ll & \frac{\log^2 x}{x} \cdot \frac{x^{1/2} \log x }{ \log (x^{1/2} \log x )} 
 \ll \frac{\log^2 x}{x^{1/2}} 
\end{eqnarray*}
by the Prime Number Theorem.

Thus if
$s_{n} > \epsilon / n$ for $ \gg x^{1/2} \log^3 x$ natural numbers $x < n \leq 2x$,
then $$ \sum_{x< n < 2x} s_{n} \gg \frac{1}{x} (x^{1/2} \log^3 x) = \frac{\log^3 x}{x^{1/2}}  ,$$ 
a contradiction. 
By dyadic decomposition, $s_{n} = o(1/n)$
holds for $\ll x^{1/2} \log^3 x = o(x/ \log x)$ numbers 
$n \leq x$.
\end{proof}

Recall $g_{n} = \lim_{x \to \infty} g_{n} (x)$, where 
$$
g_{n}(x) = \sum_{\substack{q \ \text{prime} \\ q \leq x \\ q \equiv 1 \mod n }} \frac{1}{q} \ 
- \sum_{\substack{q \ \text{prime} \\ q \leq x \\ q \equiv -1 \mod n }} \frac{1}{q}.
$$
Then by 
Proposition \ref{grprop1nn} we have $f_n = g_n + o(1/n)$ for all but 
$o(x/\log x)$ numbers 
$n \leq x$.

\begin{prop}
\label{grprop2nn}
Assume the Elliott-Halberstam Conjecture is true, and fix $ \delta > 0$. For a constant $C \geq 3$, the equation 
$g_{n} - g_{n} (n^{1 + \delta}) = o(1/n) $ holds for all but  $\ll x / \log^C x$  
natural numbers $n \leq x$. 
\end{prop}

\begin{proof}

Set $S(t, x) = \sum_{x <n< 2x} | \pi (t, n, 1) - \pi (t, n, -1) |$; then the Elliott-Halberstam Conjecture gives 
\begin{eqnarray*}
\sum_{x< n < 2x} |g_{n} - g_{n} (n^{1 + \delta})| 
& \ll & \left[ \frac{S(t, x)}{t} \right]_{x^{1+ \delta}}^{\infty} 
         + \int_{x^{1+ \delta}}^{\infty} \frac{S(t, x)}{t^2} dt \\
& \ll & \frac{1}{\log^5 x} + \int_{x^{1+ \delta}}^{\infty} \frac{dt}{t \log^5 t}
 \ll  \frac{1}{\log^{A-1} x}. 
\end{eqnarray*}
Take $A \geq 3$ in the Elliott-Halberstam Conjecture. 

If the inequality 
$$
\sum_{x< n < 2x} |g_{n} - g_{n} (n^{1 + \delta})| > \frac{\epsilon}{n}$$ 
holds for $ \gg x / \log^{A-2} x$ numbers $x < n \leq 2x$, then 
$$
\sum_{x< n < 2x} |g_{n} - g_{n} (n^{1 + \delta})|  \gg
\frac{1}{x} \cdot \frac{x}{\log^{A-2} x} = \frac{1}{\log^{A-2} x},
$$
a contradiction. 
The result follows by dyadic decomposition.
\end{proof}

\begin{cor}
\label{grcor1nn}
Assume the Elliott-Halberstam Conjecture. Then for any $\delta >0$ and $C\geq 3$, 
$$f_{n} = g_{n}(n^{1 + \delta }) + o(1/n) $$ 
for all but 
$o(x / \log x)$ numbers $n \leq x$.  
\end{cor}

\begin{proof}
This follows from Propositions \ref{grprop1nn} and \ref{grprop2nn}.
\end{proof}

We will also need the following result.
\begin{lem}[\cite{hr}, Theorem 5.7]
Let $g$ be a natural number, and let $a_i$, $b_i (1 = 1, \dots , g)$ be integers satisfying
$$
E:= \prod_{i=1}^g a_i \prod_{1 \leq r < s \leq g} (a_r b_s - a_s b_r) \neq 0.
$$ 
Let $\rho (p)$ denote the number of solutions of 
$$
\prod_{i=1}^g (a_i n +b_i) \equiv 0 \mod p,
$$
and suppose that $\rho (p) < p$ for all $p$.
Let $y$ and $x$ be real numbers satisfying $1 < y \leq x$.
Then
$$| \{n: x-y <n \leq x, a_i n +b_i \ \text{prime for} \ i=1, \dots, g \}|$$
$$\leq 2^g g! \prod_p \left( 1 - \frac{\rho (p) - 1}{p-1} \right) \left( 1- \frac{1}{p}  \right)^{-g+1} \frac{y}{\log ^g y} $$
$$ \times  \left\{  1 +  O \left( \frac{\log \log 3y + \log \log 3|E|}{\log y} \right) \right\},$$
where the implied constant depends at most on $g$.
\label{halberstam 5.7} 
\end{lem}

We are now ready to prove Theorem \ref{nnwkccv}.

\begin{proof}

We are left to deal with the finite sum
$$
g_n (n^{1+\delta}) = \sum_{q < n^{1+\delta}} \frac{c_n (q) }{q} 
= \sum_{q <(1+ \delta) n \log^2 n} \frac{c_n (q) }{ q}
+\sum_{ (1+ \delta) n \log^2 n <q < n^{1+\delta}} \frac{c_n (q) }{q}.
$$

For $ X < n < 2X$, and
$(1 + \delta )^A n \log ^A n < x < n^{1+ \delta }$,
the conditions of Lemma \ref{hooley} are satisfied, and 
\begin{eqnarray*}
\sum_{ (1+ \delta )^A n \log^A n <q < n^{1+\delta}} \frac{c_n (q) }{ q} 
& \ll & \int_{(1 + \delta )^A n \log ^A n}^{n^{1+ \delta }} \frac{t}{\phi (n) \log n} \ \frac{dt}{t^2} \\ 
& \ll & \frac{1}{\phi (n) \log n} [\delta \log n - A \log (1+ \delta ) - A \log n ] \\
& \ll & \frac{\delta }{\phi (n)}.
\end{eqnarray*}

Now we put bounds on the sum
$$
\sum_{q <(1+ \delta)^A n \log^A n} \frac{c_n (q) }{ q}  
.$$

For the range $(1+ \delta) n \log^2 n < q  <(1+ \delta)^A n \log^A n$, we use the Brun-Titchmarsh Theorem. 
\begin{eqnarray*} 
\sum_{\substack{ (1+ \delta) n \log^2 n < q \\ <(1+ \delta)^A n \log^A n}} \frac{c_n (q) }{ q}  
& \ll & \int _{(1+ \delta) n \log^2 n}^{(1+ \delta)^A n \log^A n} \frac{dt}{\phi (n) t \log (t/n)} 
\ll  \frac{1}{\phi (n)}.
\end{eqnarray*}

For $\epsilon n \log n < q  < (1+ \delta) n \log^2 n$, 
\begin{eqnarray*}
\sum_{ \substack{ \epsilon n \log n < q \\  < (1+ \delta) n \log^2 n}} \frac{c_n (q) }{q}
& \ll & \sum_{\substack{ \epsilon n < t < (1+ \delta ) \log^2 n, \\ nt \pm 1 \ \text{prime} }} \frac{1}{nt} .
\end{eqnarray*}
On average,
\begin{eqnarray*}
\sum_{x<n<2x} \ \ \sum_{\substack{ \epsilon \log n < t < (1+ \delta )\log^2 n, \\ nt \pm 1 \ \text{prime} }} \frac{1}{nt} 
& \ll & \sum_{\epsilon \log x < t < (1+ \delta ) \log^2 x } \frac{1}{tx} 
	\sum_{\substack{ x<n<2x, \\ nt \pm 1 \ \text{prime}}} 1 \\
& \leq & \sum_{\epsilon \log x < t < (1+ \delta ) \log^2 x } \frac{1}{tx} \frac{x}{\phi (t) \log x} \\ 
& \ll & \frac{\log \log x} {\log x}.
\end{eqnarray*}
Therefore the number of $ n \leq x $ such that 
$$\sum_{\substack{ \epsilon \log n < t < \log^2 n, \\ nt \pm 1 \ \text{prime} }} \frac{1}{nt}  > \frac{\epsilon}{n} $$ 
is 
$\ll \frac{x \log \log x}{ \log x} = o (x).$
That is $|g_n(\epsilon n \log n) - g_n((1+ \delta) n \log^2 n) | = o(n) $ for all but $o(x)$ natural numbers 
$n \leq x$.

Finally, for the range $q < \epsilon n \log n$, we will use Lemma \ref{halberstam 5.7}. 
Fix some $t <  \epsilon \log x$. Then the number of $ n < x$ such that $nt +1$ or $nt -1$ is a prime is 
$\leq  x / \log x +O (\frac{x \log \log x} {\log ^2 x})$, and so 
$$
\sum_{t <\epsilon \log x} \# \{n < x : \ nt \pm 1  \ \text{prime} \} \leq \epsilon x + O \left( \frac{x \log \log x}{\log x}  \right).
$$

So the number of $n < x$ for which 
$$\sum_{q < \epsilon n \log n} \frac{c_n (q) }{q} \neq 0$$ is $o(x)$, and we may assume
$g_n(\epsilon n \log n) = 0 $ for almost all $n$,  
and the proof of Theorem \ref{nnwkccv} is complete. 
\end{proof}

\section{Conditional Disproof of Generalised Kummer's Conjecture}

Recall that 
$$
h_n^- = G(n) \exp \left( \frac{\phi (n)}{2} f_n \right)
$$
and that the Generalized Kummer's Conjecture predicts that 
$h_n^- \sim G(n)$ as $n \to \infty$. 
In this section, we will prove Theorem 
\ref{gkcf}, that 
the Elliott-Halberstam Conjecture implies that 
the Generalized Kummer's Conjecture fails 
for infinitely many natural numbers $n$. 
Here we wish to show that 
$ f_n = o (1/\phi (n))$
fails for infinitely many $n$.

By Corollary \ref{grcor1nn}, the 
Elliott-Halberstam Conjecture implies  
that for any $\delta >0$ and $C\geq 3$, 
$$f_n = g_n(n^{1+\delta}) + o (1/n)$$
for all but 
$o(x / \log x)$ numbers $n \leq x$.
We wish to find bounds on $g_n(n^{1 + \delta } ) - \frac{1}{n+1}$ and so estimate the contribution of 
the primes $q$ of the form $n+1$. 
\begin{lem}
\label{gkcflem1}
Fix $\lambda > 0$ and $\epsilon > 0$. There exists some $\delta >0$ such that for all
sufficiently large values of $x$, there are 
$\leq \frac{\lambda x}{\log x} $ natural  numbers $n \leq x$ such that $n+1$ is prime and 
$$
\left|g_n(n^{1 + \delta } ) - \frac{1}{n+1} \right| \geq \frac{\epsilon }{2n} .  
$$
\end{lem}

\begin{proof}
Define 
\begin{equation*}
N_k^{\pm} (x) =  | \left\{ x<n \leq 2x : n+1 \text{ and } kn \pm 1 \ \text{both prime} \right\} | .
\end{equation*}
Then  for $k \geq 2$,
\begin{equation}
\label{N_k^+}
N_k^+ (x) \ll \left( \prod_{p | k(k-1)} \frac{p}{p-1} \right) \frac{x}{\log^2 x}
\end{equation}
and 
\begin{equation}
\label{N_k^-}
N_k^- (x) \ll \left( \prod_{p | k(k+1)} \frac{p}{p-1} \right) \frac{x}{\log^2 x}
\end{equation}
by Theorem \ref{mp 2.3}.

Thus
\begin{eqnarray*}
\sum_{\substack{x<n \leq 2x \\ n+1 \ \text{prime}}} n \left|g_n(n^{1 + \delta } ) - \frac{1}{n+1} \right| 
& \ll &  \sum_{\substack{ x<n \leq 2x \\ n+1 \ \text{prime} }} 
         \sum_{\substack{ q \equiv \pm 1 \mod n \\ q \leq n^{1+ \delta} \\ q \ \text{prime}, \ q \neq n+1  }} 
         \frac{n}{q \mp 1} \\
& \ll & \sum_{k=2}^{(2x)^{\delta}} \frac{N_k^+}{k}  + \sum_{k=2}^{(2x)^{\delta}} \frac{N_k^-}{k} \\
& \ll & \frac{x}{\log^2 x} \sum_{k=2}^{(2x)^{\delta}} \left( \frac{1}{k} \prod_{p | k(k-1)} \frac{p}{p-1} \right) 
 + \frac{x}{\log^2 x} \sum_{k=2}^{(2x)^{\delta}} \left( \frac{1}{k} \prod_{p | k(k+1)} \frac{p}{p-1} \right) \\
& \ll & \frac{x}{\log^2 x} \cdot \delta \log x = \frac{\delta x}{\log x} 
\end{eqnarray*}
by Lemmas \ref{mp 2.3} and \ref{mp 2.4}.

Hence there exists a constant $c_1 >0$ such that 
$$
\sum_{\substack {x<n \leq 2x \\ n+1 \ \text{prime}}} \left|g_n(n^{1 + \delta } ) - \frac{1}{n+1} \right|  
\leq \frac{1}{n} c_1 \delta \frac{x}{\log x} = \frac{\epsilon }{2n} \cdot \frac{\lambda x}{\log x}
$$
for $\delta = \frac{\epsilon \lambda}{2 c_1}.$ By dyadic decomposition, there exists 
$\leq \frac{\lambda x}{\log x} $ numbers $n \leq x$ such that $n+1$ is prime and 
$$
\left|g_n(n^{1 + \delta } ) - \frac{1}{n+1} \right| \geq \frac{\epsilon }{2n} ,  
$$
as required.
\end{proof}

\begin{prop}
\label{gkcfprop1}
Suppose the Elliott-Halberstam Conjecture is true. Then for a fixed $\epsilon > 0$ there exist $\gg x / \log x$ 
natural numbers $ n \leq x$ such that $f_n = (1 \pm \epsilon) \frac{1}{n}$.  
\end{prop}

\begin{proof}
By the Prime Number Theorem, there are $\geq c_2 x / \log x $ numbers $n \leq x$ such that $n+1$ is prime. 
Fix $\epsilon >0$ and let $\lambda =c_2 / 2$ in Lemma \ref{gkcflem1}. Then 
there exist $\geq c_2 x / 2 \log x$ numbers $n \leq x$ such that $n+1$ is prime and 
$$
\left|g_n(n^{1 + \delta } ) - \frac{1}{n+1} \right| < \frac{\epsilon }{2n}.  
$$

By Corollary \ref{grcor1nn}, there exist $\geq c_2 x / 2 \log x$ numbers $n \leq x$ such that $n+1$ is prime and 
$$
\left| f_n - \frac{1}{n} \right| < \frac{\epsilon}{n}. 
$$
\end{proof}

We have shown that with the assumption of the Elliott-Halberstam Conjecture, 
$$
\frac{h_n^-}{G(n)} = \exp \left( \frac{\phi (n)}{2} f_n  \right)
= \exp \left( \frac{\phi (n)}{4n} (1 \pm \epsilon) \right)
$$
for $\gg \frac{x}{\log x}$ numbers $n \leq x$. We now wish to prove that for infinitely 
many of these $n$, $\frac{\phi (n)}{n}$ is bounded away from $0$. This must be verified, of course, because  
$\liminf_{n \to \infty} \frac{\phi (n)}{n} =0$.

\begin{lem}
\label{murtylem1}
$$
\sum_{\substack{ n \leq x \\ 2n+1 \ \text{prime} }} \frac{\phi (n)}{n}
\sim \frac{c x}{\log x}, \ \ \ \text { as $x \to \infty $, }
$$
where $c= \frac{3}{2}\prod_{p \ \text{odd}} \left( 1 - \frac{1}{p(p-1)} \right) \neq 0$. 
\end{lem}

\begin{proof}
Since 
$$\frac{\phi(n)}{n} = \sum_{d|n} \frac{\mu (d)}{d} , $$
we have 
$$
\sum_{\substack{ n \leq x \\ 2n+1 \ \text{prime} }} \frac{\phi (n)}{n} \ \ =
 \sum_{\substack{ n \leq x \\ 2n+1 \ \text{prime} }} \sum_{d|n} \frac{\mu (d)}{d} \ \ = \ \
\sum_{d \leq x} \frac{\mu (d)}{d}  \sum_{\substack{ t \leq x/d \\ 2dt+1 \ \text{prime} }} 1 .
$$
The inner sum is $\pi (2x+1, 2d, 1)$. 

Thus we need to evaluate 
$$
\sum_{d \leq x} \frac{\mu (d)}{d} \pi(2x+1, 2d, 1) = 
\sum_{d < \log^A x} \frac{\mu (d)}{d} \pi(2x+1, 2d, 1) + \sum_{\log^A x \leq d \leq x} \frac{\mu (d)}{d} \pi(2x+1, 2d, 1).
$$
For the second sum, we have the estimate
$$
\sum_{\log^A x \leq d \leq x} \frac{\mu (d)}{d} \pi(2x+1, 2d, 1)
\ll \sum_{d \geq \log^A x} \frac{x}{d^2} \ll \frac{x}{\log^A x} .
$$
For the first sum, we use the Siegel-Walfisz Theorem \ref{sw} to get 
$$
\sum_{d < \log^A x} \frac{\mu (d)}{d} \pi(2x+1, 2d, 1) \ll \sum_{d \leq \log^A x} \frac{\mu (d)}{d \phi (d)} \li 2x + O \left(  \frac{x}{\log^A x}  \right).
$$
The first term is 
$$
(\li 2x) \left( \sum_{d=1}^{\infty} \frac{\mu (d) }{ d \phi (2d) } + O \left( \frac{1}{\log^{A-1} x} \right) \right).
$$
Since 
\begin{eqnarray*}
\sum_{d=1}^{\infty} \frac{\mu (d)}{d \phi (d) } 
& = & \sum_{\substack{d= 2d_1 \\ d_1 \ \text{odd} }} \frac{\mu (d)}{d \phi (2d) } 
      +  \sum_{d \ \text{odd} } \frac{\mu (d)}{d \phi (2d) }  \\
& = & \sum_{d \ \text{odd} } \frac{\mu (2d)}{2d 2\phi (d) } +\sum_{d \ \text{odd} } \frac{\mu (d)}{d \phi (d) } 
\  = \  - \frac{1}{4} \sum_{d \ \text{odd} } \frac{\mu (d)}{d \phi (d) } +\sum_{d \ \text{odd} } \frac{\mu (d)}{d \phi (d) } \\
& = & \frac{3}{4} \sum_{d \ \text{odd} } \frac{\mu (d)}{d \phi (2d) }  
\ =  \ \frac{3}{4} \sum_{d \ \text{odd}} \mu (d) \prod_{p|d} \left( \frac{p}{p-1} \right) \frac{1}{d^2} \\
& = & \frac{3}{4} \prod_{p \ \text{odd}} \left( 1 - \frac{1}{p(p-1)}   \right)
\end{eqnarray*}
and 
$$
\li 2x \sim \frac{2x}{\log x},
$$
this completes the proof of the lemma. 
\end{proof}

\begin{lem}
\label{murtylem2}
Let $\delta > 0$ and $x \geq 1$. Then
\begin{equation}
\label{crazy sum} 
\sum_{\substack{ n \leq x \\ n+1 \ \text{prime} \\ kn \pm 1 \ \text{both prime for some} 
\\ 2 \leq k \leq \delta \log x }}
\frac{\phi (n)}{n}
\ll \frac{\delta x}{\log x},
\end{equation}
where the implied constant is absolute.
\end{lem}

\begin{proof}
Since $\frac{\phi(n)}{n} \leq 1$, we have that the sum in (\ref{crazy sum}) is bounded by 
$$
\sum_{2 \leq k \leq \delta \log x} N_k^{\pm} (x)
$$
where 
\begin{eqnarray*}
N_k^{\pm} (x) & = & | \left\{ n \leq x : n +1 \ \text{ and }
kn \pm 1 \ \text{both prime} \right\}|. 
\end{eqnarray*}
By dyadic decomposition and equations (\ref{N_k^+}) and (\ref{N_k^-}), 
\begin{equation*}
N_k^+ (x) \ll \left( \prod_{p | k(k-1)} \frac{p}{p-1} \right) \frac{x}{\log^2 x}
\end{equation*}
and 
\begin{equation*}
N_k^- (x) \ll \left( \prod_{p | k(k+1)} \frac{p}{p-1} \right) \frac{x}{\log^2 x}.
\end{equation*}
Inserting this estimate into the sum, we get the sum is 
$
\ll \frac{\delta x}{\log x}, 
$
as claimed. 
\end{proof}

We now prove Theorem \ref{gkcf}.
\begin{proof}
By Proposition \ref{gkcfprop1}, the Elliott-Halberstam Conjecture implies that   
$$
\frac{h_n^-}{G(n)} = \exp \left( \frac{\phi (n)}{2} f_n  \right)
= \exp \left( \frac{\phi (n)}{2n} (1 \pm \epsilon) \right)
$$
for $\gg \frac{x}{\log x}$ numbers $n \leq x$. By Lemmas \ref{murtylem1} and \ref{murtylem2}, 
for a sufficiently small $\delta$ we have  
$$
\sum_{\substack{ \frac{x}{2} \leq n \leq x  \\ n+1 \ \text{prime} \\ \text{neither } kn \pm 1 \ \text{is prime for any} 
\\ 2 \leq k \leq \delta \log x } } 
\frac{\phi (n)}{n} \gg \frac{x}{\log x}.
$$
From this we deduce that there are infinitely many $n$ such that 
$$
\frac{h_n^-}{G(n)} \geq \exp (\eta)
$$
for some fixed $\eta > 0$. Thus if the Elliott-Halberstam Conjecture is true, then 
the Generalized Kummer's Conjecture (\ref{gkc}) fails for infinitely many natural numbers.
\end{proof}

\end{document}